# QUANTILE PYRAMIDS FOR BAYESIAN NONPARAMETRICS

By Nils Lid Hjort and Stephen G. Walker[1]

*University of Oslo and University of Kent*

Pólya trees fix partitions and use random probabilities in order to construct random probability measures. With quantile pyramids we instead fix probabilities and use random partitions. For nonparametric Bayesian inference we use a prior which supports piecewise linear quantile functions, based on the need to work with a finite set of partitions, yet we show that the limiting version of the prior exists. We also discuss and investigate an alternative model based on the so-called substitute likelihood. Both approaches factorize in a convenient way leading to relatively straightforward analysis via MCMC, since analytic summaries of posterior distributions are too complicated. We give conditions securing the existence of an absolute continuous quantile process, and discuss consistency and approximate normality for the sequence of posterior distributions. Illustrations are included.

**1. Introduction and summary.** Constructing manageable classes of random probability measures is at the heart of nonparametric Bayesian methodology. Recent surveys of Bayesian nonparametric methods, including description of several such classes of random distributions, have been given in Walker et al. (1999) and Hjort (2003). The aim of the present article is to introduce and investigate one more such class, namely that of quantile pyramids.

One attempt to construct a random probability measure on $[0,1]$ is via so-called Pólya trees. This relies on the idea of a fixed binary tree partition of $[0,1]$ and a strategy for allocating random mass to these partitions. The original and clearest exposition is provided by Ferguson (1974). More recent work on Pólya trees has been done by Lavine (1992, 1994). Inference is attractively simple since, given an independent and identically distributed

Received December 2006; revised August 2007.

[1]Supported by an Advanced Research Fellowship from the Engineering and Physical Sciences Research Council, UK.

*AMS 2000 subject classifications.* 60G35, 62F15.

*Key words and phrases.* Consistency, Dirichlet process, nonparametric Bayes, Bernshteĭn–von Mises theorem, quantile pyramids, random quantiles.

set of observations, the posterior is also a Pólya tree and the update is straightforward. A drawback to Pólya trees, and perhaps the main reason why they have not seen much application within the Bayesian nonparametric literature, is that an arbitrary partition tree of $[0,1]$ needs to be specified. There is no obvious selection criterion, though on $[0,1]$ the dyadic intervals are the natural choice. No partition is "right," however, and two different partitions produce two different answers. No satisfactory solution to this problem can be anticipated.

The fundamental idea of the Pólya tree is a fixed partition and random mass. We turn this around and instead use the idea of fixed mass and random partitions. The arbitrariness is now lost as the quantiles form a nonarbitrary partition of mass. For a distribution with cumulative function $F$ on $\mathcal{R}$, the quantile function is

$$(1) \qquad Q(y) = F^{-1}(y) = \inf\{t : F(t) \geq y\} \qquad \text{for } 0 < y < 1.$$

Our program is to construct random probability distributions $F$ via their quantile functions $Q$, using $F(x) = \sup\{y : Q(y) \leq x\}$. Specifically, the first random partition at $Q(\frac{1}{2})$ corresponds to the median and the fixed mass of $\frac{1}{2}$ is allocated in equal measure to $[0, Q(\frac{1}{2}))$ and $[Q(\frac{1}{2}), 1]$. The random partitions at $Q(\frac{1}{4})$ and $Q(\frac{3}{4})$ on the second level determine the quartiles and the fixed mass of $\frac{1}{4}$ is allocated to the relevant intervals. At stage three we draw the octiles $Q(\frac{1}{8})$, $Q(\frac{3}{8})$, $Q(\frac{5}{8})$, $Q(\frac{7}{8})$. In general, at level $m$, we draw quantiles $Q(j/2^m)$ for $j = 1, 3, \ldots, 2^m - 1$. Even though more general probabilistic constructions could be envisaged, we focus on those pyramidal schemes where $Q(j/2^m)$ for $j = 1, 3, \ldots, 2^m - 1$ are drawn independently, conditional on the values generated at level $m-1$ above, with $Q(j/2^m) \in [Q((j-1)/2^m), Q((j+1)/2^m))$.

A mild disadvantage of our quantile trees is that the prior to posterior computation is not analytically tractable, or at any rate less so than for Pólya trees. However, with the recent advent of simulation based inference the need for clear-cut conjugacy and analytically tractable posteriors is no longer critical. We shall rely on simulation strategies to collect samples from the posterior distribution. Therefore, we do not see the lack of analytical tractability as a problem and we have removed the need to specify an arbitrary partition. The allocation of the fixed quantile masses to the random partitions is the obvious choice, since they are instantly recognizable and interpretable.

While nonparametric priors are typically difficult to manipulate, in the sense that the incorporation of real qualitative prior information is nontrivial, we believe the contrary is true for quantile pyramids. The significance of quantiles is well understood and hence assigning a prior to the median, quartiles, etc. should be relatively straightforward. There are instances in the literature suggesting that more of statistics, from modeling to analysis and interpretation, should be carried out using quantiles; see, for example, Parzen (2004, 1979).

The layout of the paper is as follows. In Section 2 we introduce the quantile pyramid process on $[0,1]$. In particular we discuss issues of existence and continuity. That the pyramid schemes have a large nonparametric support is demonstrated in Section 3. In Section 4 we consider, in particular, the Beta quantile pyramid.

In Section 5 we proceed with Bayesian inference associated with quantile pyramids. First we use the quantile pyramids to construct a prior on the space of piecewise linear quantile functions. We undertake exact posterior inference for such priors for any finite level of the pyramid. We also consider a multinomial type pseudo-likelihood function for the quantiles, and investigate the implied pseudo-posterior distribution of the parameters of a quantile pyramid. The pseudo-likelihood function in question is a natural generalization of a suggestion of Jeffreys (1967), Section 4.4, concerning the median parameter, and is sometimes called the substitution likelihood; cf. Lavine (1995) and Dunson and Taylor (2005).

Then in Section 6 we work out the structure of the posterior quantile pyramids, given a set of independent data points. It is shown that the likelihood functions factorize in precisely the same way as the quantile pyramid priors, leading to simplifications of the posteriors. We demonstrate how to obtain summaries from the posterior quantile pyramid via MCMC algorithms. In Sections 7 and 8 results about the large-sample behavior of the posterior distributions are obtained; in particular Bernshteĭn–von Mises type theorems are proved under natural conditions. Finally, in Section 9 we provide a brief discussion with concluding remarks.

**2. Quantile pyramid processes.** This section considers ways of assigning a probability distribution to the full quantile process, and investigates conditions under which it is absolutely continuous. For simplicity of presentation we work on the unit interval, and consider therefore processes $\{Q(y): 0 \leq y \leq 1\}$ with $Q(0) = 0$ and $Q(1) = 1$. Such a $Q$ process is linked to a cumulative distribution function $F$ via (1). Note that $Q$ is the left-continuous inverse of the right-continuous $F$, and that $Q(y) \leq x$ if and only if $y \leq F(x)$. This somewhat nontrivial equivalence is valid also for cases where $F$ has jumps; see, for example; Shorack and Wellner (1986), Chapter 1.1.

2.1. *General pyramid quantile processes.* Consider a quantile process down to level $m$, involving random quantiles $Q(j/2^m)$ for $j = 1, \ldots, 2^m - 1$. We say that $Q$ is a pyramid quantile process down to this level if these $2^m - 1$ quantiles have been generated by successive conditionally independent mechanisms down to level $m$. More specifically, this corresponds to having the median $Q(\frac{1}{2})$ drawn from some density $\pi_{1,1}$ on $[0,1]$; then the two quartiles $Q(\frac{j}{4})$ for $j=1,3$ drawn independently from two densities $\pi_{2,1}, \pi_{2,3}$ concentrated on respectively $[0, Q(\frac{1}{2})]$ and $[Q(\frac{1}{2}), 1]$; then the four remaining octiles $Q(\frac{j}{8})$ for $j=1,3,5,7$ independent from four level-three distributions $\pi_{3,j}$ confined to the appropriate intervals $[0, Q(\frac{1}{4})]$, $[Q(\frac{1}{4}), Q(\frac{1}{2})]$, $[Q(\frac{1}{2}), Q(\frac{3}{4})]$, $[Q(\frac{3}{4}), 1]$; and so on. The simultaneous density of the $2^m - 1$ quantiles can therefore be represented as

$$\pi_{1,1}\Big(Q\Big(\frac{1}{2}\Big)\Big) \pi_{2,1}\Big(Q\Big(\frac{1}{4}\Big) \,\Big|\, Q\Big(\frac{1}{2}\Big)\Big) \pi_{2,3}\Big(Q\Big(\frac{3}{4}\Big) \,\Big|\, Q\Big(\frac{1}{2}\Big)\Big)$$

(2)
$$\times \prod_{j=1,3,5,7} \pi_{3,j}\Big(Q\Big(\frac{j}{8}\Big) \,\Big|\, \text{parents}\Big) \cdots$$

$$\times \prod_{j=1,3,5,\ldots,2^m-1} \pi_{m,j}\Big(Q\Big(\frac{j}{2^m}\Big) \,\Big|\, \text{parents}\Big)$$

where the parents of $Q(j/2^m)$ are $Q((j\pm 1)/2^m)$, both of whom were created in the previous generation.

2.2. *Existence and absolute continuity.* We now examine the quantile pyramid building process in some more detail, where variables at level $m$ are generated after those of level $m-1$. At this level,

(3)  $Q_m(j/2^m) = Q_{m-1}((j-1)/2^m)(1 - V_{m,j}) + Q_{m-1}((j+1)/2^m) V_{m,j}$

for $j = 1, 3, 5, \ldots, 2^m - 1$, in terms of independent variables $V_{m,j}$'s at work at level $m$ of the process. Note that variables $V_{m,j}$ at level $m$ are allowed to depend on previous generations' $V_{m',j'}$ for $m' \leq m-1$. We define $Q_m$ on the full unit interval by linear interpolation outside the $j/2^m$ points, and with $Q_m(0) = 0$, $Q_m(1) = 1$. Under various sets of conditions there will be a well-defined process $Q$ to which $Q_m$ converges in distribution, in the space $D_L[0,1]$ of left-continuous functions with right-hand limits on the unit interval, equipped with the Skorohod topology; see Billingsley (1968), Chapter 4, for definitions. We shall outline two arguments that can be used to establish existence of and convergence to $Q$.

The first line of arguments uses martingales. For simplicity of presentation assume now that the $V_{m,j}$'s of (3) all have mean $\frac{1}{2}$; more general results follow with additional efforts. Then, for each $y$, $Q_m(y)$ forms a martingale sequence with respect to the history up to and including the parents, and $\mathrm{E} Q_m(y) = y$. Hence there is a limit $Q(y)$ to which $Q_m(y)$ converges with probability 1. Clearly, the limit $Q$ is nondecreasing, with $Q(0) = 0$ and $Q(1) = 1$, that is, a random quantile function. We note that $Q$ in some cases might not be continuous almost surely. Such martingale arguments are pursued further in Propositions 2.2 and 2.3 below.

The second line of arguments involves tightness and is less immediate, but provides more information, in particular continuity, when the criterion we develop now applies. The $V_{m,j}$ variables of (3) are now allowed to be fully general, without the mean $\frac{1}{2}$ constraint.

PROPOSITION 2.1. *Assume that*
$$\Delta_m = \max_{j \le 2^m} \{Q_m(j/2^m) - Q_m((j-1)/2^m)\} \to_p 0.$$
*Then there is a well-defined random continuous quantile process $Q$ to which $Q_m$ converges, in the space $C[0,1]$ of continuous functions on the unit interval, equipped with the uniform topology.*

PROOF. The crux is that the condition given implies tightness of the $\{Q_m\}$ sequence in the $C[0,1]$ space, as we demonstrate in the next paragraph. Given the tightness, Prokhorov's theorem secures the existence of a subsequence converging in distribution to a limit process $Q$, which also must be continuous; see Billingsley (1968), Chapter 2. The values of this limit process at dyadic points are identical to those of $Q_m$. By denseness of dyadic points it follows that also other subsequences must have the same limit, proving that $Q$ is the limit process of $Q_m$.

To prove tightness it suffices by the theory of Billingsley (1968), Chapter 2, to show that for each positive $\varepsilon$ and $\varepsilon'$, there is a $\delta$ such that
$$\Pr\{\omega(Q_m, \delta) \ge \varepsilon\} \le \varepsilon' \qquad \text{for all large } m,$$
where $\omega(Q_m, \delta)$ is the maximum of all $\delta$-increments $Q_m(y') - Q_m(y)$ with $y' - y \le \delta$. Now let $\delta = (\frac{1}{2})^m$. For such $y$ and $y'$, find dyadic neighbors with $i/2^m \le y \le y' \le j/2^m$, where $j - i$ is at most 2. Hence, for all $m' \ge m$,
$$Q_{m'}(y') - Q_{m'}(y) \le Q_{m'}(j/2^m) - Q_{m'}(i/2^m) \le 2\Delta_m,$$
using that $Q'_m$ is equal to $Q_m$ at all $j/2^m$ points. It follows that
$$\Pr\{\omega(Q_{m'}, (\tfrac{1}{2})^m) \ge \varepsilon\} \le \Pr\{2\Delta_m \ge \varepsilon\} \qquad \text{for all } m' \ge m,$$
proving tightness under the $\Delta_m \to_p 0$ condition. □

We now show that the condition of Proposition 2.1 is very easily fulfilled. Assume, for example, that all the $V_{m,j}$ variables are independent, and that all $EV_{m,j}^2$ and $E(1-V_{m,j})^2$ are bounded by some $c < \frac{1}{2}$—if, for example, the $V_{m,j}$ is Beta with parameters $(\frac{1}{2}a_m, \frac{1}{2}a_m)$, then the condition holds provided only that the $a_m$'s stay away from zero. Note that $Q_m(j/2^m) - Q_m((j-1)/2^m)$ may be expressed as a generic product $V_1^{\varepsilon_1} \cdots V_m^{\varepsilon_m}$, with $V_1$ from the first generation, $V_2$ from the second, etc.; and where $\varepsilon_i$ is 0 or 1, writing on this occasion $V^0 = 1 - V$ and $V^1 = V$. See Figure 1. See also (5) below. This leads to

$$\Pr\{\Delta_m \geq \varepsilon\} \leq \sum_{j \leq 2^m} \Pr\{Q_m(j/2^m) - Q_m((j-1)/2^m) \geq \varepsilon\}$$
$$\leq \sum_{j \leq 2^m} (1/\varepsilon^2) \mathrm{E}(V_1^{\varepsilon_1})^2 \cdots (V_m^{\varepsilon_m})^2 \leq (1/\varepsilon^2)(2c)^m,$$

showing that Proposition 2.1 applies. Similarly, if all $V_{m,j}$'s have their means inside $[0.293, 0.707]$ and if their variances go to zero, then the $\Delta_m \to_p 0$ condition holds, implying again a continuous quantile limit process $Q(y)$. Examples of such random $Q(y)$ curves are presented in Figure 1.

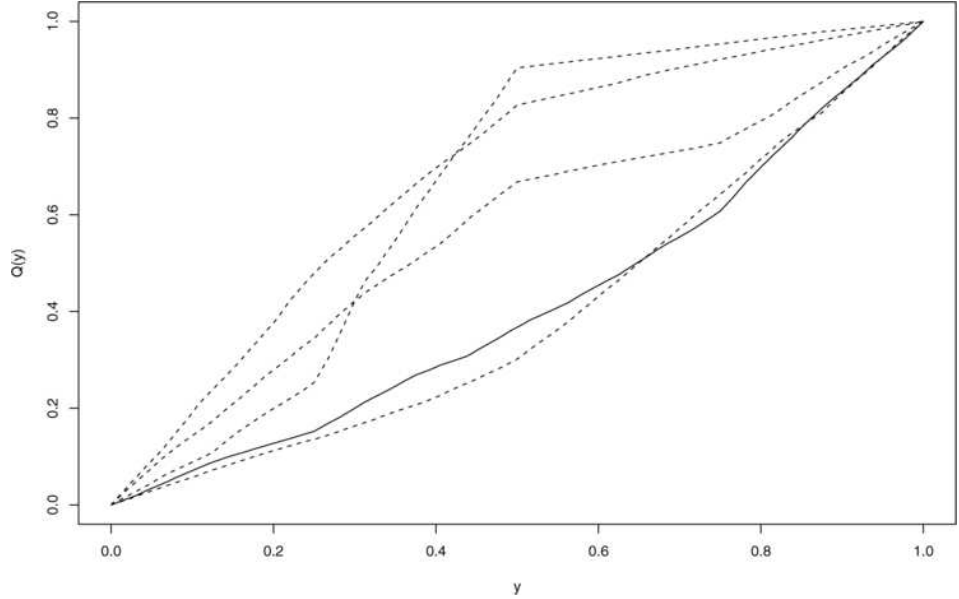

FIG. 1.  *Random $Q(y)$ curves generated from the same quantile pyramid process, using independent $V_{m,j}$'s drawn from Beta distributions $(\frac{1}{2}a_m, \frac{1}{2}a_m)$, with $a_m = cm^3$, for $c = 2.5$. This corresponds to the Beta quantile pyramid discussed in Section 4, and produces absolutely continuous $Q(y)$.*

The behavior of the $Q$ process depends crucially on aspects of the $V_{m,j}$ variables. Now we focus on conditions securing smoothness of the $Q$, and for which we must demand more than Proposition 2.1. Consider therefore the derivative of $Q_m$ at level $m$, which exists outside the $j/2^m$ points;

$$q_m(y) = \{Q_m(j/2^m) - Q_m((j-1)/2^m)\}/(\tfrac{1}{2})^m$$
(4)
$$\text{on } ((j-1)/2^m, j/2^m).$$

We wish to establish conditions under which this quantile density function converges to a random function which may be represented as the derivative of $Q$. For illustration, take $m = 3$, where we may write

$$
(5) \quad q_3(y) = \begin{cases} 8V_{1,1}V_{2,1}V_{3,1}, & \text{for } y \in (0, 1/8), \\ 8V_{1,1}V_{2,1}(1 - V_{3,1}), & \text{for } y \in (1/8, 2/8), \\ 8V_{1,1}(1 - V_{2,1})V_{3,3}, & \text{for } y \in (2/8, 3/8), \\ 8V_{1,1}(1 - V_{2,1})(1 - V_{3,3}), & \text{for } y \in (3/8, 4/8), \\ 8(1 - V_{1,1})V_{2,3}V_{3,5}, & \text{for } y \in (4/8, 5/8), \\ 8(1 - V_{1,1})V_{2,3}(1 - V_{3,5}), & \text{for } y \in (5/8, 6/8), \\ 8(1 - V_{1,1})(1 - V_{2,3})V_{3,7}, & \text{for } y \in (6/8, 7/8), \\ 8(1 - V_{1,1})(1 - V_{2,3})(1 - V_{3,7}), & \text{for } y \in (7/8, 1). \end{cases}
$$

We shall see that increased tightness of the $V_{m,j}$'s around $\frac{1}{2}$ as $m$ grows is the key to a well-behaved limit of $q_m(y)$. In fact we shall now state and prove two results securing existence of an absolutely continuous limiting quantile process $Q$, the first for the symmetric case where the $V_{m,j}$'s have mean $\frac{1}{2}$ and the second for the nonsymmetric case.

PROPOSITION 2.2. *Assume that the variables $V_{m,j}$ of (3) involved at level $m$ are such that $\mathrm{E}(V_{m,j} \mid \mathcal{F}_{m-1}) = \frac{1}{2}$ and $\mathrm{Var}(V_{m,j} \mid \mathcal{F}_{m-1}) \leq \sigma_m^2$ for each $j$, with $\sum_{m=1}^{\infty} \sigma_m^2$ finite, where $\mathcal{F}_{m-1}$ represents all previous $V_{m',j'}$ with $m' \leq m - 1$. Then with probability 1 there is a function $q(y)$ which is the a.e. limit of $q_m(y)$, and for which $Q(y) = \int_0^y q(u)\, du$ for $0 \leq y \leq 1$.*

PROOF. As seen above, each increment $Q_m(j/2^m) - Q_m((j-1)/2^m)$ at level $m$ may be represented as a product $V_1^{\varepsilon_1} \cdots V_m^{\varepsilon_m}$, where $\varepsilon_i$ is 0 or 1, and $V^0 = 1 - V$, $V^1 = V$. Hence $q_m(y)$ may be presented as $Z_m = W_1 \cdots W_m$ with $W_j = 2V_j^{\varepsilon_j}$ having mean 1. Thus the martingale convergence theorem applies and leads to the existence of a limit $q(y)$, regardless of the variances.

The finiteness of $\sum_{m=1}^{\infty} \sigma_m^2$ is, however, needed in order to secure that $Q(y)$ is the integral of $q(y)$. The variance of $Z_m$ above is $\mathrm{E}W_1^2 \cdots W_m^2 - 1$, which via conditional expectations and $\mathrm{E}(W_m^2 \mid \mathcal{F}_{m-1}) \leq 1 + 4\sigma_m^2$ is seen to be bounded by $\prod_{r=1}^m (1 + 4\sigma_r^2) - 1$, which again is bounded by the constant $\exp(4 \sum_{m=1}^{\infty} \sigma_m^2) - 1$. In particular, $\int_0^1 \mathrm{Var}\, q_m(y)\, dy$ is bounded as $m$ grows. The required statement follows from the corollary of Kraft (1964). The point to note is that as far as independence of the $V_{m,j}$'s is concerned, it is their conditional independence given $\mathcal{F}_{m-1}$ which actually matters, and with this the theorem and corollary of Kraft (1964) still hold, since $Q$ behaves as a distribution function with probability 1. □

The above quantile processes with $\mathrm{E}(V_{m,j} \mid \mathcal{F}_{m-1}) = \frac{1}{2}$ all have $\mathrm{E}Q(y) = y$ for $y \in (0, 1)$, that is, are centered at the uniform quantile function. In Bayesian practice one needs to be able to center priors at given positions, that is, to adjust the prior to match a given prior guess distribution, say $Q_{\text{null}}(y) = F_{\text{null}}^{-1}(y)$. This takes nonsymmetric $V_{m,j}$'s. One needs in fact

$$\text{E}V_{m,j} = \frac{Q_{\text{null}}(j/2^m) - Q_{\text{null}}((j-1)/2^m)}{Q_{\text{null}}((j+1)/2^m) - Q_{\text{null}}((j-1)/2^m)}$$

(6)
$$\text{for } j = 1, 3, \ldots, 2^m - 1$$

at level $m$, as is seen from representation (3), in order to achieve $\text{E}Q_m(y) = Q_{\text{null}}(y)$. If the distribution with $Q_{\text{null}}$ as quantile function has a density with two derivatives, an approximation to the mean above is

$$\tfrac{1}{2} - \tfrac{1}{4}\{Q''_{\text{null}}(y)/Q'_{\text{null}}(y)\}/2^m,$$

at $y = j/2^m$. The following proposition demonstrates that such nonsymmetric setups also give absolutely continuous quantile pyramids, provided the variances become small enough.

We also point out another option for achieving a similar aim, via a simple transformation, namely through $Q(y) = Q_{\text{null}}(Q_{\text{unif}}(y))$, where $Q_{\text{unif}}(y)$ is a quantile process centered at the uniform distribution, using symmetric $V_{m,j}$'s. The median of this random $Q(y)$ is equal to $Q_{\text{null}}(y)$. For example, $Q(y) = \mu + \sigma\Phi^{-1}(Q_{\text{unif}}(y))$ defines a quantile process with median value function equal to the quantile function of a normal $(\mu, \sigma^2)$, with $\Phi$ denoting the cumulative distribution function of a standard normal.

PROPOSITION 2.3. *Assume that the $V_{m,j}$'s are all independent, and write $\text{E}V_{m,j} = \tfrac{1}{2} + \delta_{m,j}$ and the unconditional variance $\text{Var}\,V_{m,j} = \sigma_{m,j}^2$. Assume further that $|\delta_{m,j}| \le \delta_m$ and $\sigma_{m,j} \le \sigma_m$ for all $j$ at level $m$, where $\sum_{m=1}^{\infty}\sigma_m^2$ and $\sum_{m=1}^{\infty}\delta_m$ are both finite. Then again there is a.s. convergence of $q_m(y)$ to $q(y)$, and $Q(y)$ is the integral of $q(y)$.*

PROOF. As in the previous proof we may represent $q_m(y)$ as $Z_m = W_1 \cdots W_m$, with $W_k = 2V_k^{\varepsilon_k}$, again writing generically $V^1 = V$ and $V^0 = 1 - V$. Existence of a limit for $Z_m$ is not as automatic as in the previous symmetric case, since its mean differs from 1 and the martingale convergence theorem cannot be directly applied. Consider, however, $Z_m^* = W_1 \cdots W_m/(\xi_1 \cdots \xi_m)$, where $\xi_1 \cdots \xi_m$ is the mean of $Z_m$. Then the martingale theorem applies to $Z_m^*$, which therefore has a well-defined limit $Z^*$. But the sequence of products of means $\xi_1 \cdots \xi_m$ is seen to converge, basically since the conditions imposed imply that $\prod_{r=m+1}^{m'} \xi_r$ must converge to 1 when $m$ and $m' > m$ go to infinity. Next,

$$\text{E}Z_m^2 \le \prod_{r=1}^{m}\{(1+2\delta_r)^2 + 4\sigma_r^2\} = \prod_{r=1}^{m}(1+2\delta_r)^2 \prod_{r=1}^{m}\left\{1 + \frac{4\sigma_r^2}{(1+2\delta_r)^2}\right\}.$$

Using the $1+x \leq \exp(x)$ inequality and some further analysis, one sees that the sequence of variances is also bounded, like the sequence of means. That $Q(y)$ is the integral of the limiting quantile density function $q(y)$ under the boundedness of $\operatorname{Var} q_m(y)$ condition follows as for the previous proposition, again via techniques from the proof of Kraft's (1964) theorem. $\square$

**3. Large pyramidal support.** Assume that sufficient conditions are in place for $Q$, and hence also $F = Q^{-1}$, to be absolutely continuous with respect to the Lebesgue measure on $[0,1]$; cf. the proposition above. Assume the same to be true for $Q_0$, which we shall refer to as the true quantile function. Then $Q_0$ admits the density $f_0$ on $[0,1]$ with corresponding quantile function $Q_0(y) = F_0^{-1}(y)$ and quantile density $q_0(y) = 1/f_0(Q_0(y))$. Now let $\Pi$ be the probability measure governing the $q = \lim q_m$, and consider the following conditions:

(A) For all $\varepsilon > 0$, $\Pi\{q : \int q \log(q/q_0) \, du < \varepsilon\} > 0$.
(B) For all $\delta > 0$ there exists an $\varepsilon > 0$ such that

$$\int \log \frac{q_0(\tau_\varepsilon(u))}{q_0(u)} \, du < \delta$$

for any $\tau_\varepsilon(u)$ for which $\max_u |\tau_\varepsilon(u) - u| < \varepsilon$ and $\tau_\varepsilon(u) \in [0,1]$.
(C) The density $f_0$ is bounded by some $K < \infty$.

PROPOSITION 3.1. *When conditions* (A)–(C) *hold, each Kullback–Leibler neighborhood* $\{f : \int f_0 \log(f_0/f) \, dx < \varepsilon\}$ *around the fixed* $f_0$ *has positive $\Pi$-probability.*

PROOF. We first show that condition (B) implies condition (B0), which is that for all $\delta > 0$, there exists an $\varepsilon > 0$ such that $\int \{q_0(u)/q_0(\tau_\varepsilon(u))\} \, du > 1 - \delta$ for any $\tau_\varepsilon(u)$ for which $\max_u |\tau_\varepsilon(u) - u| < \varepsilon$. For any $\theta > 0$ there exists an $\varepsilon > 0$ such that

$$\int \log \frac{q_0(u)}{q_0(\tau_\varepsilon(u))} \, du > -\theta.$$

Now for any positive random variable $Z$, $\log \operatorname{E} Z > \operatorname{E} \log Z$ and so

$$\log \int \frac{q_0(u)}{q_0(\tau_\varepsilon(u))} \, du > -\theta$$

and hence

$$\int \frac{q_0(u)}{q_0(\tau_\varepsilon(u))} \, du > \exp(-\theta).$$

Clearly, for any $\delta > 0$ there exists a $\theta > 0$ such that $\exp(-\theta) > 1 - \delta$ and so the claim is proven.

By using the transform $x = Q_0(u)$, the Kullback–Leibler distance $\int f_0 \log(f_0/f)\,dx$ may be expressed as

$$\int \log \frac{q(\tau(u))}{q_0(u)}\,du = \int \log \frac{q(\tau(u))}{q_0(\tau(u))}\,du + \int \log \frac{q_0(\tau(u))}{q_0(u)}\,du,$$

where $\tau(u) = F(Q_0(u))$. We will now deal with these two terms separately.

For the first term, use the transform $x = \tau(u)$ to give

$$\int q(x) \log\{q(x)/q_0(x)\} f_0(Q(x))\,dx.$$

The aim is to show that for any $\delta > 0$, the prior puts positive mass on

$$\int q(x) \log\{q(x)/q_0(x)\} f_0(Q(x))\,dx \leq K \int q(x) \log\{q(x)/q_0(x)\}\,dx + \delta$$

and hence, using condition (A), the prior puts positive mass on

$$\int q(x) \log\{q(x)/q_0(x)\} f_0(Q(x))\,dx < \delta$$

for any $\delta > 0$. Now

$$\int q(x) \log\{q(x)/q_0(x)\}\{K - f_0(Q(x))\}\,dx \geq \int f_0(Q(x)) q_0(x)\,dx - 1$$

and

$$\int f_0(Q(x)) q_0(x)\,dx = \int \frac{q_0(x)}{q_0(\lambda(x))}\,dx,$$

where $\lambda(x) = F_0(Q(x))$. Condition (A) is sufficient for the prior to put positive mass on $\{Q : \max_u |Q(u) - Q_0(u)| < \theta\}$ for any $\theta > 0$ and so from the absolute continuity of $F_0$, for any $\varepsilon > 0$ the prior puts positive mass on $\{Q : \max_u |\lambda(u) - u| < \varepsilon\}$. Condition (B0) finishes the story for the first term.

For the second term, again, condition (A) is sufficient for the prior to put positive mass on $\{Q : \max_u |Q(u) - Q_0(u)| < \theta\}$ for any $\theta > 0$. Thus, from the absolute continuity of $F$, for any $\varepsilon > 0$ the prior puts positive mass on $\{F : \max_u |\tau(u) - u| < \varepsilon\}$. Hence, using condition (B), for any $\delta > 0$ there exists an $\varepsilon > 0$ such that the prior puts positive mass on

$$\int \log \frac{q_0(\tau_\varepsilon(u))}{q_0(u)}\,du < \delta.$$

This completes the proof. □

Here we establish that condition (A) holds in all situations where the $V_{m,j}$'s have full support on $[0,1]$, expectations fixed at $\frac{1}{2}$, and with variances decreasing sufficiently fast. It will be clear from the arguments used that condition (A) continues to be in force also when the expectations deviate slightly from $\frac{1}{2}$, within the limits dictated by Proposition 2.3. The sufficiently fast variances in question are recorded in Barron, Schervish and Wasserman (1999) and amount to

$$\sum_{m=1}^{\infty} \max_j (\operatorname{Var} V_{m,j})^{1/2} \leq \sum_{m=1}^{\infty} \sigma_m < \infty.$$

With the conditional expectations of the $V_{m,j}$'s fixed at $\frac{1}{2}$ it is noted that

$$\operatorname{E}\operatorname{Var}(V_{m,j} \mid \mathcal{F}_{m-1}) = \operatorname{Var} V_{m,j}$$

and so conditioning on previous $V_{m',j}$ is permitted provided

$$\sum_{m=1}^{\infty} \max_j \{\operatorname{E}\operatorname{Var}(V_{m,j} \mid \mathcal{F}_{m-1})\}^{1/2} < \infty,$$

which is the sum of the unconditional variances, since $\operatorname{E}(V_{m,j}) = \frac{1}{2}$, in the sense that the argument of Barron, Schervish and Wasserman (1999) then continues to go through. Following Lavine (1994), we can write $D(Q, Q_0) = \int q \log(q/q_0) \, dy$ as the difference of two sums, the first being

$$(7) \qquad \sum_m \sum_\varepsilon \sum_{j \in \{0,1\}} Q(B_{\varepsilon,j}) \log \frac{Q(B_{\varepsilon,j} \mid B_\varepsilon)}{\lambda(B_{\varepsilon,j} \mid B_\varepsilon)},$$

where $\lambda$ represents the Lebesgue measure. The $\lambda$ can also be replaced with any smooth $Q_{\text{null}}$, as long as it remains a dominating measure for the $Q$ process; cf. the previous section. Here $B_\varepsilon$ is a dyadic interval and for a particular $m$, we have $\varepsilon = (\varepsilon_1, \ldots, \varepsilon_m)$ where $\varepsilon_k \in \{0, 1\}$ and $(\varepsilon, j) = (\varepsilon_1, \ldots, \varepsilon_m, j)$ for $j \in \{0, 1\}$. So, for example, $B_{0,0} = [0, 1/4]$ and $B_{0,1,1} = [3/8, 4/8]$. Now make the variances of the $V_\varepsilon$ decrease sufficiently rapidly, ensuring that

$$\sum_m \max_\varepsilon \left\{ \log \frac{V_\varepsilon}{\lambda(B_{\varepsilon,j} \mid B_\varepsilon)} \vee \log \frac{1 - V_\varepsilon}{\lambda(B_{\varepsilon,j} \mid B_\varepsilon)} \right\}$$

converges with positive probability and hence that (7) converges. The second term is $\int q \log q_0 \, dy$ which is finite if (7) is, and which can also be expressed as a sum over $m$,

$$\sum_m \sum_\varepsilon \sum_{j \in \{0,1\}} Q(B_{\varepsilon,j}) \log \frac{Q_0(B_{\varepsilon,j} \mid B_\varepsilon)}{\lambda(B_{\varepsilon,j} \mid B_\varepsilon)}.$$

The proof is then completed using the fact that the $V$'s have full support on $[0, 1]$.

We note that the result about condition (A) being satisfied is not surprising in view of the fact that the $Q$ is similar enough in structure to a Pólya tree in order for Theorem 2 of Lavine (1994) to apply.

Condition (B) is a property of $q_0$ and is a quite mild smoothness condition. If $q_0$ is continuous on $(0, 1)$, then we only need consider the integral in a neighborhood of 0 and 1.

**4. The Beta and the median-Dirichlet quantile pyramids.** We shall discuss two attractive options for the class of distributions $V_{m,j}$ of (3) that make up the core of a quantile pyramid.

4.1. *The Beta quantile pyramid.* The first option is to use independent Betas for the weights appearing in (3). An appealing choice is to let all $V_{m,j}$'s in (3) be symmetric Beta variables, with parameters, say, $(\frac{1}{2}a_m, \frac{1}{2}a_m)$ for those at level $m$. These have variance $\frac{1}{4}/(a_m+1)$. As long as $\sum_{m=1}^{\infty} 1/(a_m+1)$ is finite, the limiting quantile process $Q$ of the $Q_m$'s has a.s. a quantile density function $q(y)$. The slightly stronger condition $\sum_{m=1}^{\infty} 1/a_m^{1/2} < \infty$ secures condition (A) of Section 3. Note that the $Q$ is constructed as a Pólya tree, but importantly its accompanying distribution function $F = Q^{-1}$ is not. Interestingly, the very same condition about the Beta distribution parameters, about $\sum_{m=1}^{\infty} 1/a_m^{1/2}$ being finite, occurs in Ghosal, Ghosh and Ramamoorthi (1999b), where it is seen to imply posterior consistency of symmetrized Pólya trees.

The uncertainty of $q(y)$ around its constant mean 1 is dictated by the variances of the $V_{m,j}$'s, sometimes in complicated ways. Intriguingly, when all $V_{m,j}$'s inside the same generation $m$ have the same distribution, symmetric around $\frac{1}{2}$, we may actually find and assess the distribution of $B_m = \max_y q_m(y)$ and its limit $B = \max q(y)$ explicitly. This is due to the symmetry of the representation $2V_1^{\varepsilon_1} \cdots 2V_m^{\varepsilon_m}$ over different intervals, as displayed, for example, in (5). At each node, either $V_{m,j}$ or $1 - V_{m,j}$ is in $(\frac{1}{2}, 1)$, the other in $(0, \frac{1}{2})$. The maximum value of $q_m(y)$ takes place in that interval for which each of the $m$ components $V_j^{\varepsilon_j}$ are in $(\frac{1}{2}, 1)$. Hence

$$B_m =_d \prod_{j=1}^{m} \max\{2V_j, 2(1-V_j)\} =_d \prod_{j=1}^{m}(2U_j)$$

in terms of generic $V_1, V_2, \ldots$ from generations $1, 2, \ldots$, where $U_j$ is distributed like $V_j$ conditional on $V_j \geq \frac{1}{2}$. This distribution converges, under the finite sum of variances condition, and is easily simulated for given regimes for the distributions of $V_j$'s. For the Beta quantile pyramid, $EU_m = \xi(\frac{1}{2}a_m)$, say, where $\xi(b)$ is the mean of $V \mid \{V \geq \frac{1}{2}\}$ when $V$ is a symmetric $\text{Beta}(b, b)$. Some efforts and integration skills lead to

$$\xi(b) = \int_{1/2}^{1} v \frac{\Gamma(2b)}{\Gamma(b)^2} 2v^{b-1}(1-v)^{b-1} dv = \frac{\Gamma(2b)}{\Gamma(b)^2} \left(\frac{1}{4}\right)^b \left\{\frac{1}{b} + \frac{\Gamma(1/2)\Gamma(b)}{\Gamma(b+1/2)}\right\}.$$

Hence $EB_m = \prod_{j=1}^{m} \{2\xi(\frac{1}{2}a_j)\}$.

To assess this usefully, we note that $2(2b+1)^{1/2}(V - \frac{1}{2})$ tends to a standard normal, when $V \sim \text{Beta}(b, b)$ and $b$ goes to infinity. This means that $U = \max(V, 1-V)$ behaves like $\frac{1}{2} + \frac{1}{2}Z/(2b+1)^{1/2}$ for large $b$, where $Z$ is a standard normal conditional on being positive. This leads to $2\xi(b) \doteq 1 + (2/\pi)^{1/2}/(2b+1)^{1/2}$ for increasing $b$. We learn from this that $B = \max q(y)$ has finite mean when $\sum_{m=1}^{\infty} 1/a_m^{1/2}$ is finite. Also, if, for example, $a_m = cm^3$, then the mean of $B$ may be studied as a function of $c$, which is useful when attempting to elicit a prior process for one's quantile function. One may similarly study the distribution and expected value of $\min q(y) =_d \prod_{j=1}^{\infty} 2(1 - U_j)$.

As discussed around (6), we would often wish to center quantile processes at given null distributions, which would need nonsymmetric Beta variables in the above construction, say, employing $\text{Beta}(\frac{1}{2}a_m, \frac{1}{2}b_m)$ with appropriate $a_m, b_m$. Proposition 2.3 dictates that $a_m$ and $b_m$ need to become close to each other for growing $m$, in order for a limiting quantile density function $q(y)$ to exist.

One special case worth mention is that where all the $V_{m,j}$'s are uniform, corresponding to all $a_m = 2$, where the quantile process amounts to a natural splitting procedure: (i) the median $Q(\frac{1}{2})$ is uniform on $[0,1]$; (ii) the two extra quartiles are independent and uniform on $[0, Q(\frac{1}{2})]$ and $[Q(\frac{1}{2}), 1]$; (iii) the three extra octiles are independent and uniform on the four intervals defined by the three quartiles; and so on. This might be seen as a natural noninformative prior scheme. More generally one might study the case of $a_m = a$ constant, with the same $\text{Beta}(\frac{1}{2}a, \frac{1}{2}a)$ at work at all levels for the $V_{m,j}$. Then the $Q(y)$ process is a.s. continuous but singular, not equal to the integral of its derivative. This follows from results of Ferguson (1974), page 621; see also Dubins and Freedman (1967), based on the fact that $Q$ behaves as a distribution function with probability 1.

The $V_{m,j}$'s of the Beta quantile pyramid might employ parameters $(\frac{1}{2}a_m, \frac{1}{2}b_m)$ that depend on the previous outcomes of $V_{m',j}$ for $m' < m$, for example, in a Markovian fashion. This gives one the opportunity to modify the behavior of $Q_m$ in light of aspects of $Q_{m-1}$.

4.2. *The median-Dirichlet quantile pyramid.* Agree to say that a random variable $U$ has a median-Dirichlet distribution with parameter $a$, written $U \sim \text{MD}(a)$, if

(8) $\Pr\{U \leq x\} = H_a(x) = \Pr\{\text{Beta}(ax, a(1-x)) \geq \frac{1}{2}\} = G(\frac{1}{2}; a(1-x), ax).$

Here $G(\cdot; a, b)$ denotes the cumulative Beta distribution with parameters $(a, b)$. To motivate this definition, suppose that $F$ is a Dirichlet process with parameter $aF_{\text{uni}}$, where $F_{\text{uni}}$ is the uniform distribution on the unit interval. Then its random median $U = Q(\frac{1}{2}) = \inf\{t : F(t) \geq \frac{1}{2}\}$ does in fact have this $\text{MD}(a)$ distribution. Note that $U$ is symmetric around its center value $\frac{1}{2}$. More generally, when $F_{\text{null}}$ is any probability distribution on the line, say that $U \sim \text{MD}(aF_{\text{null}})$ if $\Pr\{U \leq x\} = H_a(F_{\text{null}}(x))$. This is the distribution of a random median from a $\mathrm{Dir}(aF_{\mathrm{null}})$ process; see also Hjort and Petrone (2006). Our emphasis in this connection is on using the MD$(a)$ distribution as a modeling tool when working with quantile pyramids.

From Proposition 2.2 we know that the degree of continuity of the limiting quantile process is governed by the sizes of the variances of the $V_{m,j}$'s. For the present case this necessitates studying the variances of the MD$(a)$ distribution (8). This may be written

$$\tau^2(a) = \int_0^1 \mathrm{Pr}_a\{U^2 \geq x\}\,dx - (\tfrac{1}{2})^2 = \int_0^1 G(\tfrac{1}{2}; ax^{1/2}, a(1-x^{1/2}))\,dx - 1/4.$$

Inspection and some analysis reveal that $\tau^2(a)$ starts with value $1/12$ for $a$ at zero and then goes down at rate $O(1/a)$ when $a$ grows. Intriguingly, $\tau^2(a) = (1/4)\rho(a)/(a+1)$ where $\rho(a)$ goes monotonically from $1/3$ up to $1$ as $a$ grows, making the MD$(a)$ distribution quite similar to the Beta$(\tfrac{1}{2}a, \tfrac{1}{2}a)$ for growing $a$. Hence remarks made earlier for the Beta quantile pyramids have clear analogues for the median-Dirichlet governed quantile pyramids; convergence of $\sum_{m=1}^\infty 1/(a_m+1)$ secures absolute continuity, for example. It may also be attractive to determine the concentration parameter $a$ at level $m$ by taking into account the results realized at level $m-1$. One such option is $V_{m,j} \mid \mathcal{F}_{m-1} \sim \mathrm{MD}(a_{m,j})$ with $a_{m,j} = b_m/\Delta(m-1,j)$, in terms of $\Delta(m-1,j) = Q_{m-1}((j+1)/2^m) - Q_{m-1}((j-1)/2^m)$. Finiteness of $\sum_{m=1}^\infty 1/(1+b_m)$ secures absolute continuity of the resulting quantile pyramid.

**5. Exact posterior and pseudo-posterior pyramids.** Let $X_1, \ldots, X_n$ be independent observations from a continuous distribution $F$ on $[0,1]$. We shall discuss ways of obtaining the posterior distribution of the quantile process.

One point of view is that $Q$ defines the cumulative distribution function $F$, after which aspects of the posterior distribution of $F$ may in principle be derived via the defining characteristics

$$\mathrm{Pr}\{F \in C, X_1 \in A_1, \ldots, X_n \in A_n\} = \mathrm{E}I\{F \in C\}F(A_1)\cdots F(A_n),$$

valid for all Borel subsets $C$ of the space of cumulative distribution functions and for all intervals $A_1, \ldots, A_n$. Then aspects of $Q$ given data may be derived using (1). For example, considering a single $Q(y)$,

$$\mathrm{Pr}\{Q(y) \leq x, X_i \in x_i \pm \varepsilon \text{ for each } i\}$$
$$= \mathrm{E}I\{Q(y) \leq x\} \prod_{i=1}^n \{Q^{-1}(x_i + \varepsilon) - Q^{-1}(x_i - \varepsilon)\},$$

which in principle should lead to the posterior distribution of $Q(y)$. This would often be a cumbersome route to follow, however, which is why we circumvent the $F$ here, attempting instead to work directly with the $Q$ process.

First a comment on how we do this is in order. We obtain the exact posterior for the prior $\Pi_m$. Such a prior generates random quantile functions by linear interpolation. We work with $\Pi_m$ in the same way that posterior Pólya trees are obtained for partitions constructed down to a finite level. What is important is that the prior is defined and exists for all $m$ and as $m \to \infty$ converges to a well-defined prior $\Pi$.

5.1. *Exact posterior inference for $\Pi_m$.* The prior, as has been mentioned generates quantile functions based on linear interpolation between random points. Then the inverse of $Q_m$, say $F_m$, is linear on each quantile interval $[q_{j-1}, q_j]$, with a constant derivative there;

$$(9) \qquad f_m(x) = F'_m(x) = \frac{1}{k} \frac{1}{q_j - q_{j-1}} \qquad \text{for } x \in (q_{j-1}, q_j),$$

for $j = 1, \ldots, k = 2^m$. Here $q_0 = 0$ and $q_k = 1$. This amounts to a "random histogram" type model, with random cell widths but fixed probabilities over these cells.

There is also another route to the (9) density, as follows. In general, for a smooth distribution $F$ with density $f$ and quantile function $Q$, the quantile density function is $q(y) = Q'(y) = 1/f(Q(y))$. Inverting this gives $f(x) = 1/q(F(x))$. In the present context this leads naturally to the level $m$ prior which generates random densities of the type

$$f_m(x) = 1/q_m(F_m(x)),$$

where $F_m(x)$ for given $x$ is the solution $y$ to the equation $Q_m(y) = x$. But this can be seen to be exactly the same as (9), due to expression (4) for $q_m$ and the linear interpolation character of $Q_m$ and $F_m$.

Under this linear interpolation prior there is a well-defined likelihood

$$(10) \qquad \bar{L}_n(q) = \prod_{j=1}^{k} \left( \frac{1}{k} \frac{1}{q_j - q_{j-1}} \right)^{N_j(q)},$$

where $N_j(q) = nF_n(q_{j-1}, q_j]$ is the number of points falling inside the $j$th quantile interval (and with $F_n$ being the empirical distribution of the data). Its behavior for growing $n$ is dictated by

$$-n^{-1} \log \bar{L}_n(q) = \sum_{j=1}^{k} n^{-1} N_j(q) \log(q_j - q_{j-1}) + \log k$$

$$\to_p \bar{\lambda}(q) = \sum_{j=1}^{k} F_0(q_{j-1}, q_j] \log \frac{q_j - q_{j-1}}{1/k}.$$

For fixed prior and growing $n$, the posterior distribution of $q = (q_1, \ldots, q_{k-1})$ will concentrate on decreasing neighborhoods around $q^0 = (q_1^0, \ldots, q_{k-1}^0)$, defined as the minimizer of $\bar{\lambda}(q)$.

5.2. *The multinomial substitute likelihood.* In the setting above, assume that a pyramid-type probability distribution is given for the $k-1$ quantiles $q_1, \ldots, q_{k-1}$, where $k = 2^m$, but we avoid any further specification of $Q$. We define the pseudo, or substitute, likelihood for the data as the multinomial probability

$$
(11) \quad \begin{aligned} L_n(q) &= \binom{n}{N_1(q), \ldots, N_k(q)} \left(\frac{1}{k}\right)^{N_1(q)} \cdots \left(\frac{1}{k}\right)^{N_k(q)} \\ &= \frac{n!}{N_1(q)! \cdots N_k(q)!} \left(\frac{1}{k}\right)^n. \end{aligned}
$$

Such a construction can be found in Jeffreys (1967), Chapter 4, for the particular case of the median, that is, for $k = 2$, who noted that it would yield a "valid uncertainty." This has been further discussed by Kalbfleisch (1978) and by Monahan and Boos (1992), who pointed out that $L_n(q)$ is not the conditional distribution of the data given any statistic, and by Lavine (1995), who showed that in any case using this substitute likelihood produces asymptotically conservative inference. The following arguments and results provide more general insight into aspects discussed in the above references, and specifically lend support to Jeffreys's claim; further discussion is offered in Section 8.

We start with Stirling's formula and find

$$
\log \Gamma(np + 1) = (np + \tfrac{1}{2})(\log n + \log p) - np \\
+ \log(2\pi)^{1/2} + (1/12)(np)^{-1} + O((np)^{-2})
$$

for growing $np$, from which we derive

$$
-n^{-1} \log L_n(q) = \sum_{j=1}^{k} \left( \widetilde{p}_j \log \widetilde{p}_j + \frac{\log \widetilde{p}_j}{2n} \right) + \log k + \frac{k-1}{2n} \log \frac{n}{2\pi} + R_n(q),
$$

where $\widetilde{p}_j = N_j(q)/n = F_n(q_{j-1}, q_j]$ is the relative proportion of points falling inside the $j$th quantile interval and

$$
R_n(q) = \frac{1}{12} \frac{1}{n^2} \left\{ \sum_{j=1}^{k} \frac{1}{\widetilde{p}_j} - 1 \right\} + \text{smaller terms}
$$

goes to zero in probability. If the data points are generated from some $F_0$, so that $F_n \to F_0$ uniformly, with probability 1, then

$$
-n^{-1} \log L_n(q) \to_p \lambda(q) = \sum_{j=1}^{k} F_0(q_{j-1}, q_j] \log \frac{F_0(q_{j-1}, q_j]}{1/k}.
$$

This is the Kullback–Leibler distance from the discrete probability distribution with point masses $F(q_{j-1}, q_j]$ for $j = 1, \ldots, k$ to the uniform one with point masses $1/k$. This lends credibility to (11) as being appropriate for the nonparametric framework, since maximizing $L_n(q)$ for large $n$ amounts to minimizing $\lambda(q)$, which happens exactly for $F_0(q_{j-1}, q_j] = 1/k$ for each $j$, that is, $F_0(q_j) = j/k$. Also, since $\sum_{j=1}^k u_j \log(u_j/k^{-1}) \doteq \frac{1}{2} k \sum_{j=1}^k (u_j - 1/k)^2$ when all the $u_j$'s are close to $1/k$, $L_n(q)$ is approximately proportional to

$$(12) \qquad L_n^*(q) = \exp\left[-\tfrac{1}{2} nk \sum_{j=1}^k \{F_n(q_{j-1}, q_j] - 1/k\}^2\right] (\widetilde{p}_1 \cdots \widetilde{p}_k)^{1/2}.$$

**6. Updating.** Consider a pyramid quantile process of the general type described in Section 2, interpreted as a prior process for an unknown quantile function $Q_m(y)$ in a nonparametric Bayesian setup. This section describes how we may update $Q_m$ after having observed a sample $x_1, \ldots, x_n$. Following Section 5 there is the exact likelihood and the pseudo-likelihood; two related but different ways of handling the updating. We show that for both versions, the pyramid structure is retained, leading to certain simplifications for the posterior and pseudo-posterior quantile distributions.

6.1. *Updating the linear interpolation prior.* Consider a quantile process described down to level $m$, involving a $Q_m$ defined in terms of the $q_j = Q(j/k)$ quantiles for $j = 1, \ldots, k-1$, where $k = 2^m$. There is a prior $\Pi_m(q)$ for $q = (q_1, \ldots, q_{k-1})$ of the type (2). With the likelihood (10), the exact posterior is given by

$$(13) \qquad \Pi_m(q \mid \text{data}) \propto \Pi_m(q) \bar{L}_n(q).$$

We now demonstrate that the likelihood factorizes in pyramidal fashion.

The basic step involves the following quantity. Let $M_n(a, b) = n F_n(a, b]$ count the number of data points having fallen inside $(a, b]$, and study

$$\bar{\kappa}_n(q; a, b) = \left(\frac{1}{2} \frac{1}{(q-a)/(b-a)}\right)^{M_n(a,q)} \left(\frac{1}{2} \frac{1}{(b-q)/(b-a)}\right)^{M_n(q,b)}$$
$$\text{for } q \in (a, b),$$

where $M_n(a, b)$ counts the number of data points falling in $(a, b]$. To exemplify, we find for $m = 2$ that

$$\bar{L}_n(q_1, q_2, q_3) = \bar{\kappa}_n(q_2; 0, 1) \bar{\kappa}_n(q_1; 0, q_2) \bar{\kappa}_n(q_3; q_2, 1).$$

Similarly, for $m = 3$,

$$\bar{\kappa}_n\left(Q\left(\tfrac{4}{8}\right); 0, 1\right) \prod_{j=2,6} \bar{\kappa}_n\left(Q\left(\tfrac{j}{8}\right); Q\left(\tfrac{j-2}{8}\right), Q\left(\tfrac{j+2}{8}\right)\right)$$
$$\times \prod_{j=1,3,5,7} \bar{\kappa}_n\left(Q\left(\tfrac{j}{8}\right); Q\left(\tfrac{j-1}{8}\right), Q\left(\tfrac{j+1}{8}\right)\right).$$

The general formula involving $\bar{\kappa}_n(Q(j/2^m); \text{parents})$ follows and becomes

$$\bar{\kappa}_n\Big(Q\Big(\frac{1}{2}\Big); 0, 1\Big) \bar{\kappa}_n\Big(Q\Big(\frac{1}{4}\Big); 0, Q\Big(\frac{1}{2}\Big)\Big) \bar{\kappa}_n\Big(Q\Big(\frac{3}{4}\Big); Q\Big(\frac{1}{2}\Big), 1\Big)$$

$$\times \prod_{j=1,3,5,7} \bar{\kappa}_n\Big(Q\Big(\frac{j}{8}\Big); Q\Big(\frac{j-1}{8}\Big) Q\Big(\frac{j+1}{8}\Big)\Big)$$

$$\times \cdots \prod_{j=1,3,5,\ldots,2^m-1} \bar{\kappa}_n\Big(Q\Big(\frac{j}{2^m}\Big); \text{parents}\Big).$$

Verifying that this is identical to $\bar{L}_n(q)$ of (10), with $q_j = Q(j/2^m)$, is a matter of algebra and book-keeping. This leads to an expression for the posterior distribution:

$$\pi_{1,1}\Big(Q\Big(\frac{1}{2}\Big)\Big) \bar{\kappa}_n\Big(Q\Big(\frac{1}{2}\Big); 0, 1\Big)$$

(14)
$$\times \prod_{j=1,3} \pi_{2,1}\Big(Q\Big(\frac{j}{4}\Big) \,\Big|\, \text{parents}\Big) \bar{\kappa}_n\Big(Q\Big(\frac{j}{4}\Big); \text{parents}\Big)$$

$$\times \prod_{j=1,3,5,7} \pi_{3,1}\Big(Q\Big(\frac{j}{8}\Big) \,\Big|\, \text{parents}\Big) \bar{\kappa}_n\Big(Q\Big(\frac{j}{8}\Big); \text{parents}\Big) \cdots.$$

This part provides details of the Metropolis–Hastings algorithm for the linear interpolation process. The posterior density for $q = (q_1, \ldots, q_{k-1})$ is given by

$$\Pi_m(q \mid \text{data}) \propto \Pi_m(q) \prod_{j=1}^{k} \Big(\frac{1}{q_j - q_{j-1}}\Big)^{N_j(q)} \qquad \text{where } k = 2^m.$$

A Metropolis–Hastings algorithm proceeds by taking a proposal $q'$ for $q$, which we do by changing one component at a time; that is, we take $q'_j$ uniform on $(q_{j-1}, q_{j+1})$ and $q'_l = q_l$ for $l \neq j$. Consequently, the accept–reject ratio for the algorithm is

$$\min\Big\{1, \frac{(q_j - q_{j-1})^{N_j(q)}(q_{j+1} - q_j)^{N_{j+1}(q)} \Pi_m(q')}{(q'_j - q'_{j-1})^{N_j(q')}(q'_{j+1} - q'_j)^{N_{j+1}(q')} \Pi_m(q)}\Big\}.$$

This is in principle a straightforward algorithm to implement.

There is of course a broadly flexible class of priors to use when it comes to the choice of $\Pi_m(q)$, via (2). For illustration, take all the $V_{m,j}$'s of representation (3) to be independent, with the $V_{m',j}$'s at level $m'$ coming from the same density $g_{m'}$, as with the Beta quantile pyramids. Then, at level $m = 5$, $\Pi_5(q_1, \ldots, q_{31})$ may be written

$$g_1(q_{16}) \prod_{j \in S_2} g_2\Big(\frac{q_j - q_{j-8}}{q_{j+8} - q_{j-8}}\Big) \frac{1}{q_{j+8} - q_{j-8}}$$

$$
\begin{aligned}
&\times \prod_{j\in S_3} g_3\left(\frac{q_j - q_{j-4}}{q_{j+4} - q_{j-4}}\right)\frac{1}{q_{j+4} - q_{j-4}} \\
&\times \prod_{j\in S_4} g_4\left(\frac{q_j - q_{j-2}}{q_{j+2} - q_{j-2}}\right)\frac{1}{q_{j+2} - q_{j-2}} \\
&\times \prod_{j\in S_5} g_5\left(\frac{q_j - q_{j-1}}{q_{j+1} - q_{j-1}}\right)\frac{1}{q_{j+1} - q_{j-1}}
\end{aligned}
\tag{15}
$$

on the set where $0 < q_1 < \cdots < q_{31} < 1$, in which $S_2 = \{8, 24\}$, $S_3 = \{4, 12, 20, 28\}$, $S_4 = \{2, 6, 10, 14, 18, 22, 26, 30\}$, and $S_5 = \{1, 3, 5, \ldots, 31\}$.

6.2. *Updating with the multinomial substitute likelihood.* For this approach we use $L_n(q)$ instead of $\bar{L}_n(q)$ in (13), and consider

$$
\kappa_n(q; a, b) = \binom{M_n(a,b)}{M_n(a,q), M_n(q,b)} (\tfrac{1}{2})^{M_n(a,b)} \quad \text{for } q \in (a,b). \tag{16}
$$

This is also the symmetric binomial probability that $M_n(a, q)$ of the points, among the $M_n(a, b)$, will fall in the $(a, q]$ interval. Note that $M_n(a, q)$ and $M_n(q, b)$ depend in a somewhat cumbersome form on the $q$ argument. We shall see, via algebraic manipulations of the multinomial likelihood (11), that it also factories into various contributions, of the type (16). To exemplify, for the case $m = 2$ one finds that $L_n(q_1, q_2, q_3)$ equals

$$
\binom{n}{N_1(q), N_2(q), N_3(q), N_4(q)} (\tfrac{1}{4})^n = \kappa_n(q_2; 0, 1)\kappa_n(q_1; 0, q_2)\kappa_n(q_3; q_2, 1).
$$

Similarly, for $m = 3$, we find

$$
L_n(q) = \kappa_n(q_4; 0, 1)\kappa_n(q_2; 0, q_4)\kappa_n(q_6; q_4, 1) \prod_{j=1,3,5,7} \kappa_n(q_j; q_{j-1}, q_{j+1})
$$

in terms of the octiles vector $(q_1, \ldots, q_7)$, and so on. The general formula follows as for the previous case and verification is again a matter of algebra and book-keeping. This leads to an expression for the pseudo-posterior distribution of the same structure as (14).

The best way of sampling from the pseudo-posterior distribution of the vector $(q_1, \ldots, q_{k-1})$ appears to be via a Metropolis–Hastings type algorithm, as follows. A proposal for $q$ is taken to be $q'$ given by $q'_j$ uniform on $(q'_{j-1}, q'_{j+1})$ with $q'_l = q_l$ for $l \neq j$. For an iteration, we sweep through all the $j$'s in turn. Consequently, the accept–reject ratio for the algorithm is given by

$$
\min\left\{1, \frac{N_j(q)!N_{j+1}(q)!\Pi_m(q')}{N_j(q')!N_{j+1}(q')!\Pi_m(q)}\right\}.
$$

For general discussion of aspects of the Metropolis–Hastings type algorithms, see, for example, Tierney (1994).

6.3. *Illustrations.* A number of simulations were undertaken. Firstly, for a moderate sample size of $n = 100$ and with $m = 32$, the true quantile function, from which the data were simulated, was taken to be $Q(y) = y^2$. The Gibbs sampler was run for 5000 iterations and all the samples were used in constructing the Bayes estimate of $Q(y)$. Figure 2 is the Bayes estimate of $Q$ using the substitute likelihood. The bold line denotes the estimate and the dotted line the true quantile function. Figure 3 corresponds to the Bayes estimate based on the linear interpolation process. Again, the bold line is the estimate and the dotted line is the true quantile function. The prior used in both cases is the uniform for the quantile interpolators $V_{m,j}$'s, that is, as in (15) with the uniform for $g_1, g_2, g_3, g_4, g_5$. Note that the accompanying joint density for $(q_1, \ldots, q_{31})$, for this "uniform stick-breaking prior," is not flat in $q$-space.

**7. Bayesian consistency.** In this section we provide results related to Bayesian consistency and asymptotic proximity of the approaches/models used in Sections 5 and 6. We go further in Section 8, reaching large-sample approximation results of the Bernshteĭn–von Mises theorem variety.

Subject to regularity conditions on $f_0$ [see conditions (B) and (C) in Section 3], the prior can be arranged so that it puts positive mass on all Kullback–Leibler neighborhoods of $f_0$; see Proposition 3.1. With this it is well known that the posterior distributions accumulate in weak neighborhoods of $f_0$. That is, $\Pi_n(A) = \Pi(A \mid X^n) \to 1$ with $f_0$-probability 1 where $X^n = (X_1, \ldots, X_n)$ and $A$ is any weak neighborhood of $f_0$.

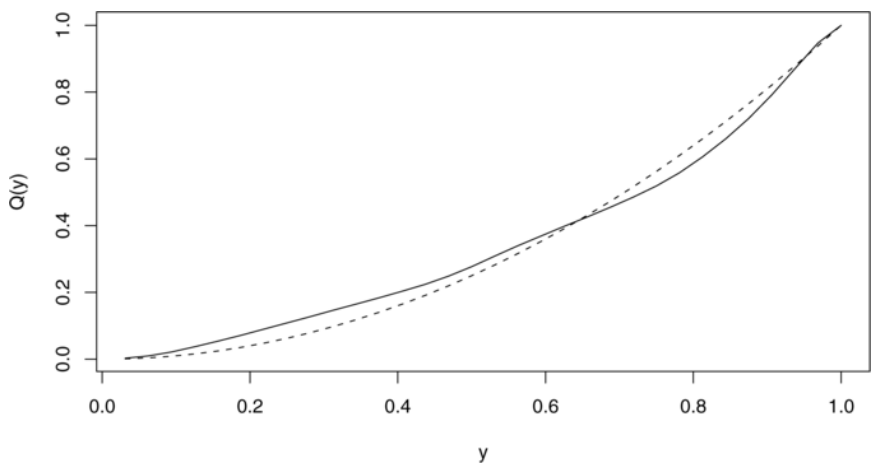

Fig. 2. *With $n = 100$ data points simulated from the distribution with quantile function $Q_0(y) = y^2$, the figure displays the Bayes estimate of $Q$ using the substitute likelihood. The bold line denotes the estimate and the dotted line the true quantile function.*

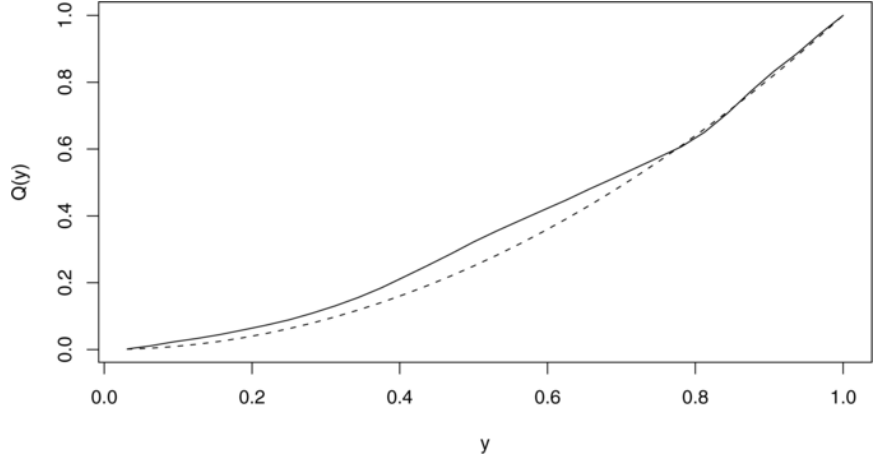

FIG. 3. *With $n = 100$ data points simulated from the distribution with quantile function $Q_0(y) = y^2$, the figure displays the Bayes estimate of $Q$ using the linear interpolation likelihood. The bold line denotes the estimate and the dotted line the true quantile function.*

Nevertheless, inference will be based on $\Pi_m$ which will typically be sample size dependent; the more samples the larger $m$ will be taken. So we will undertake consistency issues assuming that $m_n$ increases as $n$ increases. This is what an experimenter would do and so we establish consistency for such a procedure. So, consider now the prior $\Pi_n^* = \Pi_{m_n}$ which generates $f_m(x)$ defined in (9), with resolution level $m = m_n$, now allowed to increase slowly with sample size, with consequent $k = k_n = 2^{m_n}$ cells.

PROPOSITION 7.1. *Let independent observations $X_1, X_2, \ldots$ be generated from a density $f_0$, that is, inside the Kullback–Leibler support of $\Pi$, the limit of quantile pyramids $\Pi_m$, assuming conditions* (A)–(C) *of Section 3 hold. Let furthermore $\Pi_n^* = \Pi_{m_n}$ be the quantile prior stemming from construction (9) for $f_m$, with*

$$k_n \to \infty \quad \text{and} \quad k_n/n \to 0. \tag{17}$$

*Then the sequence of posterior distributions is Hellinger consistent at $f_0$.*

PROOF. For a finite $m$ the number of Hellinger balls required to fill up the space of densities generated by $\Pi_n^*$ is finite. Call this number by $N_n$. In order to achieve consistency with this sample size dependent prior, along with the support condition, as in Section 3, we require that

$$\sum_{n=1}^{\infty} \exp(-nc) N_n < \infty \quad \text{for each positive } c.$$

This result can be found in Walker (2003) based on previous findings from Ghosal, Ghosh and Ramamoorthi (1999a). The $c$ is related to the size of the balls and hence the need for it to be arbitrarily small. Let us fix the size of the balls to be $\delta > 0$. An observation is that if $|q_j - q_j^*| < \varepsilon$ for all $j = 1, \ldots, k_n - 1$, then the Hellinger distance between the corresponding densities $f$ and $f^*$ will be bounded by $\delta$ for small enough $\varepsilon$. Now split the unit interval into $K = [1/\varepsilon]$ equal parts all of size $\varepsilon$. Then clearly

$$N_n(\varepsilon) \leq K^{k_n} \qquad \text{with } k_n = 2^{m_n}.$$

So we require

$$\sum_{n=1}^{\infty} \exp\{-n(c - n^{-1} k_n \log K)\} < \infty$$

for all $c > 0$ and $K$ which happens precisely under the (17) condition. □

A refinement is possible here by allowing the size of the balls to also depend on $n$.

**8. Approximations to the pseudo-posterior distribution.** Here we examine various natural approximations to the pseudo-posterior quantile pyramids, and reach so-called Bernshteĭn–von Mises theorems under natural conditions.

For parametric models a classic large-sample result about the maximum likelihood estimator $\widehat{\theta}$ for the parameter $\theta$ is that $\sqrt{n}(\widehat{\theta} - \theta_0)$ tends to $\mathrm{N}(0, J(\theta_0)^{-1})$, with $J(\theta_0)$ the information matrix at the true parameter value $\theta_0$. A Bayesian mirror result to this is that under mild conditions on the model and the prior used for $\theta$, the posterior distribution of $\sqrt{n}(\theta - \widehat{\theta})$ will a.s. have the same limit distribution. This also implies that the Bayes estimator $\mathrm{E}(\theta \mid \text{data})$ and the maximum likelihood estimator become $\sqrt{n}$-equivalent for large $n$. Such results can be traced back to Bernshteĭn (1917) and von Mises (1931) and are often collectively referred to as Bernshteĭn–von Mises theorems. The importance of these results is partly that an easy-to-use approximation can be used in applied Bayesian statistics, in cases where the precise prior-to-posterior calculations are complicated, but lies also in revealing that data appropriately wash out the prior as the data information level increases; different priors will lead to approximately the same inference, and this inference will also agree to the first order of magnitude with that of classic frequentist approaches. In Bayesian nonparametrics such results are not to be taken for granted [see, e.g., Freedman (1999) and Hjort (2003) for counterexamples], and are also typically harder to prove if they hold at all; see, for example, Ghosal (2000), about exponential families with a growing number of parameters, and Kim and Lee (2004) and De Blasi and Hjort (2007), concerned with semiparametric event history models with Beta process type priors.

Let again $f_0$ be the true density underlying independent data $X_1,\ldots,X_n$ on $[0,1]$, with cumulative and quantile distribution functions $F_0$ and $Q_0$. Two results from classic empirical process theory are that

$$A_n(t) = \sqrt{n}\{F_n(t) - F_0(t)\} \to_d W^0(F_0(t)),$$
$$B_n(y) = \sqrt{n}\{F_n^{-1}(y) - Q_0(y)\} \to_d q_0(y)W^0(y) = W^0(y)/f_0(Q_0(y)),$$

where $W^0$ is a Brownian bridge, that is, a zero-mean normal process with covariance function $y_1(1-y_2)$ for $y_1 \le y_2$. The first convergence takes place in the space $D_R[0,1]$ of right-continuous functions with left-hand limits on $[0,1]$ while the second holds in each space $D_L[\varepsilon, 1-\varepsilon]$ of left-continuous functions with right-hand limits on $[\varepsilon, 1-\varepsilon]$, both equipped with suitable versions of the *Skorokhod* topology. For these results see, for example, Shorack and Wellner (1986), Chapter 3.

We shall first focus on quantiles $q = (q_1,\ldots,q_{k-1})$ for a fixed number $k = 2^m$ of cells, with $q_j = Q(j/k)$. For this situation the above result for $B_n$ implies for the frequentist estimator $q_j^* = F_n^{-1}(j/k)$ that

$$\sqrt{n}(q_j^* - q_j^0) \to_d q_0(j/k)W^0(j/k) \qquad \text{for } j = 1,\ldots,k-1,$$

with $q_j^0 = Q_0(j/k)$ being the real underlying quantile. Our next result provides a Bernshteĭn–von Mises mirror result to this.

PROPOSITION 8.1. *Consider any quantile pyramid prior $\Pi_m(q)$ for the quantiles $q = (q_1,\ldots,q_{k-1})$, with the number of cells $k = 2^m$ being fixed, and let the pseudo-posterior distribution of $q$ be defined in terms of the multinomial likelihood $L_n(q)$ of (11). Then with probability 1 the pseudo-posterior distribution of $q$ is such that the vector with components*

$$C_{n,j} = \sqrt{n}(q_j - q_j^*) = \sqrt{n}\{Q(j/k) - F_n^{-1}(j/k)\}$$

*converges to that of $C_j = W^0(j/k)/f_0(q_j^0)$, for $j = 1,\ldots,k-1$.*

PROOF. Write $q_j = q_j^* + \gamma_j/\sqrt{n}$ for $j = 1,\ldots,k-1$. Then

$$F_n(q_j) = F_0(q_j) + A_n(q_j)/\sqrt{n}$$
$$= F_0(q_j^*) + \{f_0(q_j^*) + A_n(q_j^*)\}/\sqrt{n} + o_p(n^{-1/2}),$$

which implies that $\sqrt{n}\{F_n(q_j) - j/k\}$ can be written

$$\sqrt{n}\{F_0(q_j^*) - j/k\} + f_0(q_j^*)\gamma_j + A_n(q_j^*) + o_p(1) = f_0(q_j^*)\gamma_j + o_p(1).$$

Consequently, with probability 1,

$$\sqrt{n}\{F_n(q_{j-1}, q_j) - 1/k\} = f_0(q_j^0)\gamma_j - f_0(q_{j-1}^0)\gamma_{j-1} + o_p(1)$$

since $q_j^* \to q_j^0$ a.s. for $j = 1, \ldots, k$. Here we write $\gamma_0 = 0$ and $\gamma_1 = 0$. From the likelihood approximation (12), the pseudo-posterior density of $(\gamma_1, \ldots, \gamma_{k-1})$ is proportional to $\Pi_m(q^* + \gamma/\sqrt{n}) L_n(q^* + \gamma/\sqrt{n})$, where, up to further factors that vanish in importance,

$$L_n(q^* + \gamma/\sqrt{n}) \doteq \exp\left[-\tfrac{1}{2}nk \sum_{j=1}^k \{f_0(q_j^0)\gamma_j - f_0(q_{j-1}^0)\gamma_{j-1}\}^2\right]$$

(18)
$$= \exp(-\tfrac{1}{2}\phi^{\mathrm{t}} \Sigma_k^{-1} \phi),$$

say; the underlying convergence is uniform over all balls $\|\gamma\| \le c$. Here $\phi$ is the $(k-1)$-vector with components $\phi_j = f_0(q_j^0)\gamma_j - f_0(q_{j-1}^0)\gamma_{j-1}$, and

$$\Sigma_k = \begin{pmatrix} k^{-1}(1-k^{-1}) & \cdots & -k^{-2} \\ & \cdots & \\ -k^{-2} & \cdots & k^{-1}(1-k^{-1}) \end{pmatrix} \quad \text{with } \Sigma_k^{-1} = \begin{pmatrix} 2k & \cdots & k \\ & \cdots & \\ k & \cdots & 2k \end{pmatrix}.$$

This implies the statement of the proposition, in view of the multinomial structure of the covariances of a Brownian bridge. □

The above implies that the pseudo-posterior quantile process

(19) $$C_n(y) = \sqrt{n}\{Q(y) - Q^*(y)\} = \sqrt{n}\{Q(y) - F_n^{-1}(y)\}$$

is such that the pseudo-posterior distribution of $C_n(y)$ tends to $C(y) = q_0(y) \times W^0(y)$ for $y$ at positions $1/k, 2/k, \ldots, (k-1)/k$, as long as the resolution level is fixed with $k = 2^m$ cells, for any pyramid prior. It also follows, by taking expectations, that the Bayes estimator $\mathrm{E}\{Q(y) \mid \mathrm{data}\}$ becomes equivalent to the frequentist estimator $F_n^{-1}(y)$ for large $n$, under these conditions, for the $y = j/k$ positions.

REMARK. We do anticipate that there is a stronger Bernshteĭn–von Mises theorem that under conditions somewhat stronger than those of Proposition 7.1 will imply that the full process $C_n$ of (19) will converge to $C = q_0(\cdot)W^0(\cdot)$, inside each of the *Skorokhod* spaces $D_L[\varepsilon, 1-\varepsilon]$. In a technical report version of the present article we have provided details of arguments that combine to formulate the conjecture that for the Beta quantile pyramid with $V_{m,j}$'s taken as $\mathrm{Beta}(\tfrac{1}{2}a_m, \tfrac{1}{2}a_m)$, and if the density $f_0$ is bounded on $[0,1]$, then conditions

$$k_n \to \infty, \qquad k_n/\sqrt{n} \to 0, \qquad ma_m/\sqrt{n} \to 0$$

secure $C_n \to_d C$ in the described sense.

**9. Discussion and concluding remarks.** We end our article with a list of concluding comments, pertaining to various aspects of our quantile pyramid processes.

9.1. *What is the $F$ of a quantile pyramid $Q$ like?* It is in general not possible to understand the distribution function $F$ generated from $Q$ in terms of, for example, analytic expressions for means or variances. It is best understood in terms of distributing mass to random partitions and relying on existence theorems, as discussed in Section 2, in analogy to what is done if $F$ had been generated by a Pólya tree.

Doss and Gill (1992) provided a machinery for bringing weak convergence results in the $F$ domain over to the $Q = F^{-1}$ domain, via compact differentiability of the inverse functional transform. Interestingly, one may now borrow their techniques to go the other way, starting with Proposition 8.1 and the $C_n$ process of (19). The result is another Bernshteĭn–von Mises theorem, stating that the posterior distribution of $\sqrt{n}(F - F_n)$, where $F$ is the random distribution function stemming from a quantile pyramid $Q$, must tend to the Brownian bridge $W^0(F_0(\cdot))$, under mild conditions.

9.2. *Semiparametric models and quantile regression.* In Section 2 we briefly pointed to quantile processes of the type $Q(y) = \mu + \sigma \Phi^{-1}(Q_{\text{uni}}(y))$, which for given $(\mu, \sigma)$ describes a prior situated at the $N(\mu, \sigma^2)$ quantile function. By in addition having a prior on $(\mu, \sigma)$ one has a semiparametric Bayesian construction for handling an uncertain distribution about the normal. The posterior distribution can be established via a Metropolis–Hastings algorithm based around the likelihood function at level $m$ given by

$$\bar{L}_n(\mu, \sigma, q) = \frac{1}{\sigma^n} \prod_{j=1}^k \left\{ \frac{1}{k} \frac{1}{\Phi^{-1}(q_{\text{uni},j}) - \Phi^{-1}(q_{\text{uni},j-1})} \right\}^{N_j(q_{\text{uni}})},$$

where now $q_{\text{uni},j} = Q_{\text{uni}}(j/2^m)$ and $N_j(q_{\text{uni}})$ is the number of observations in

$$(\mu + \sigma \Phi^{-1}(q_{\text{uni},j-1}), \mu + \sigma \Phi^{-1}(q_{\text{uni},j})).$$

Similarly one may work with quantile regression problems, of the type $Q_i(y) = a + bx_i + \sigma \Phi^{-1}(Q_{\text{uni}}(y))$, with a prior on $(a, b, \sigma)$ independent of the $Q_{\text{uni}}$ pyramid. This would be a semiparametric construction more general in spirit than that of Kottas and Gelfand (2001), who work with the Dirichlet process.

9.3. *Dependent quantile pyramids.* There are various statistically important problems associated with dependent quantile functions, for example, in finance. This might in the present context call for constructing dependent quantile pyramids, for which there are several possibilities. A particular version is as follows, elaborating on the idea that the $V_{m,j}$'s for two pyramids can be made dependent:

$$V_{m,j} = G^{-1}(\Phi(N_{m,j}); \tfrac{1}{2}a_m, \tfrac{1}{2}a_m), \qquad V'_{m,j} = G^{-1}(\Phi(N'_{m,j}); \tfrac{1}{2}a_m, \tfrac{1}{2}a_m),$$

where $N_{m,j}$ and $N'_{m,j}$ are standard normals with correlation $\rho_m$, say, and with $G$ being the cumulative distribution function for the Beta distribution. This leads to two dependent Beta quantile pyramids $Q$ and $Q'$. More generally time series of quantile pyramids can be worked with through suitable time series models for the underlying $V_{m,j}$'s.

9.4. *Asymptotics for the linear interpolation model.* Proposition 8.1 and the anticipated process generalization described in the remark following it relate to the multinomial substitute likelihood $L_n$ of (11). Results of a somewhat different nature can be reached with the linear interpolation model based likelihood $\bar{L}_n(q)$ of (10), and these modified statements need different proofs. Importantly, for a fixed fine-ness level $m$, the two quantile likelihoods (10) and (11) are concerned with two different versions of quantiles; the first is maximized by estimators that tend to the least false quantiles $q^0 = (q_1^0, \ldots, q_k^0)$ that minimize the distance function $\bar{\lambda}(q)$ of Section 5.1, whereas the second is maximized by estimators that tend to real underlying quantiles, as explained in Section 5.2. The difference between the pseudo-quantiles and real quantiles goes to zero as the level $m$ increases, however, as is implicit in Proposition 7.1. There are, however, "cube root asymptotics" that govern the large-sample behavior of distributions associated with the linear interpolation likelihood; see Hjort (2007).

9.5. *More general quantile processes.* Our $\Pi_m(q)$ priors for quantiles $Q_m$ have been constructed in a natural pyramidal fashion, and we saw in Sections 5 and 6 that the natural updating mechanisms involved likelihoods that factorized in precisely this way. More general constructions can also be worked with, however, using techniques and results from our article. One may work with $Q(1/k), \ldots, Q((k-1)/k)$ for $k$ different from the pyramid's $2^m$, with methods of Section 6 still applying, and other constructions for building suitable $\Pi_m(q_1, \ldots, q_{k-1})$ may be contemplated, as, for example, setting $q_j$ equal to $F(j/k)$ for a random distribution function $F$ on the unit interval. Our pyramids correspond to special cases of such constructions.

9.6. *Further quantilian quantities.* In our article we have developed and discussed nonparametric Bayesian tools for analyzing a quantile distribution $Q$. This is a fundamental statistical quantity, and other quantities of importance depend naturally on $Q$. Among these are the Lorenz curve

$$L(y) = \int_0^y Q(u)\,du \Big/ \int_0^1 Q(u)\,du \qquad \text{for } 0 \le y \le 1$$

and the Gini index $G = 2\int_0^1 \{1 - L(y)\}\,dy$. Since we are able to obtain posterior samples of the full $Q$ curve, with a quantile pyramid as prior, such may be used to carry out inference for, for example, the Gini index. There are also important procedures for comparing two populations in terms of their quantile functions, including Doksum's (1974) shift function $D(x) = F_2^{-1}(F_1(x)) - x$ and Parzen's (1979) comparison distribution $\pi(y) = F_2(F_1^{-1}(y))$. Here quantile pyramids may be used as priors for $Q_1 = F_1^{-1}$ and $Q_2 = F_2^{-1}$, and Bayes analysis via posterior samples of the $D(x)$ and $\pi(y)$ curves may be performed via our methods. Hjort and Petrone (2006) give a detailed analysis of these matters for the special case of Dirichlet process priors.

**Acknowledgment.** We thank Editor Joe Eaton and two anonymous reviewers for insightful and constructive comments that helped sharpen our presentation.

## REFERENCES

Barron, A., Schervish, M. J. and Wasserman, L. (1999). The consistency of posterior distributions in nonparametric problems. *Ann. Statist.* **27** 536–561. MR1714718

Bernshteĭn, S. (1917). Teoriya Veroyatnosteĭ [*Theory of Probability*]. Gostekihzdat, Moskva.

Billingsley, P. (1968). *Convergence of Probability Measures.* Wiley, New York. MR0233396

De Blasi, P. and Hjort, N. L. (2007). Bayesian survival analysis in proportional hazard models with logistic relative risk. *Scand. J. Statist.* **34** 229–257. MR2325252

Doksum, K. A. (1974). Empirical probability plots and statistical inference for nonlinear models in the two-sample case. *Ann. Statist.* **2** 267–277. MR0356350

Doss, H. and Gill, R. D. (1992). An elementary approach to weak convergence for quantile processes, with applications to censored survival data. *J. Amer. Statist. Assoc.* **87** 869–877. MR1185204

Dubins, L. E. and Freedman, D. A. (1967). Random distribution functions. *Proc. 5th Berkeley Symp. Math. Statist. Probab.* **2** 183–214. Univ. California Press, Berkeley. MR0214109

Dunson, D. B. and Taylor, J. A. (2005). Approximate Bayesian inference for quantiles. *J. Nonparametr. Statist.* **17** 385–400. MR2129840

Ferguson, T. S. (1974). Prior distributions on spaces of probability measures. *Ann. Statist.* **2** 615–629. MR0438568

Freedman, D. A. (1999). On the Bernstein–von Mises theorem with infinite-dimensional parameter. *Ann. Statist.* **27** 1119–1140. MR1740119

Ghosal, S., Ghosh, J. K. and Ramamoorthi, R. V. (1999a). Posterior consistency of Dirichlet mixtures in density estimation. *Ann. Statist.* **27** 143–158. MR1701105

Ghosal, S., Ghosh, J. K. and Ramamoorthi, R. V. (1999b). Consistent semiparametric Bayesian inference about a location parameter. *J. Statist. Plann. Inference* **77** 181–193. MR1687955

Ghosal, S. (2000). Asymptotic normality of posterior distributions for exponential families when the number of parameters tends to infinity. *J. Multivariate Anal.* **74** 49–68. MR1790613

Hjort, N. L. (2003). Topics in nonparametric Bayesian statistics (with discussion). In *Highly Structured Stochastic Systems* (P. J. Green, N. L. Hjort and S. Richardson, eds.) 455–487. Oxford Univ. Press. MR2082419

Hjort, N. L. (2007). Model selection for cube root asymptotics. In *Report No. 50/2007 from the Mathematisches Forschungsinstitut Oberwolfach, Reassessing the Paradigms of Statistical Model-Building* 33–36.

Hjort, N. L. and Petrone, S. (2006). Nonparametric quantile inference with the Dirichlet process prior. In *Advances in Statistical Modeling and Inference: Festschrift for Kjell Doksum* (V. Nair, ed.) 463–492.

Jeffreys, H. (1967). *Theory of Probability*, 3rd ed. Clarendon, Oxford.

Kalbfleisch, J. (1978). Likelihood methods and nonparametric testing. *J. Amer. Statist. Assoc.* **83** 167–170. MR0518600

Kim, Y. and Lee, J. (2004). A Bernstein–von Mises theorem in the nonparametric right-censoring model. *Ann. Statist.* **32** 1492–1512. MR2089131

Kottas, A. and Gelfand, A. (2001). Bayesian semiparametric median regression modeling. *J. Amer. Statist. Assoc.* **96** 1458–1468. MR1946590

Kraft, C. H. (1964). A class of distribution function processes which have derivatives. *J. Appl. Probab.* **1** 385–388. MR0171296

Lavine, M. (1992). Some aspects of Polya tree distributions for statistical modelling. *Ann. Statist.* **20** 1222–1235. MR1186248

Lavine, M. (1994). More aspects of Polya tree distributions for statistical modelling. *Ann. Statist.* **22** 1161–1176. MR1311970

Lavine, M. (1995). On an approximate likelihood for quantiles. *Biometrika* **82** 220–222. MR1332852

von Mises, R. (1931). *Wahrscheinlichkeitsrechnung*. Springer, Berlin.

Monahan, J. F. and Boos, D. D. (1992). Proper likelihoods for Bayesian analysis. *Biometrika* **79** 271–278. MR1185129

Parzen, E. (1979). Nonparametric statistical data modeling (with discussion). *J. Amer. Statist. Assoc.* **74** 105–131. MR0529528

Parzen, E. (2004). Quantile probability and statistical data modeling. *Statist. Sci.* **19** 652–662. MR2185587

Shorack, G. R. and Wellner, J. A. (1986). *Empirical Processes.* Wiley, New York. MR0838963

Tierney, L. (1994). Markov chains for exploring posterior distributions (with discussion). *Ann. Statist.* **22** 1701–1762. MR1329166

Walker, S. G., Damien, P., Laud, P. W. and Smith, A. F. M. (1999). Bayesian nonparametric inference for random distributions and related functions (with discussion). *J. Roy. Statist. Soc. Ser. B* **61** 485–528. MR1707858

Walker, S. G. (2003). On sufficient conditions for Bayesian consistency. *Biometrika* **90** 482–488. MR1986664

Department of Mathematics  
University of Oslo  
P.B. 1053 Blindern N-0316  
Oslo  
Norway  
E-mail: nils@math.uio.no

Institute of Mathematics, Statistics  
and Actuarial Science  
University of Kent  
Kent CT2 7NZ  
United Kingdom  
E-mail: S.G.Walker@kent.ac.uk